\documentclass[a4paper]{amsart}

\usepackage[english]{babel}
\usepackage[utf8x]{inputenc}
\usepackage[T1]{fontenc}

\usepackage[a4paper,top=3cm,bottom=2cm,left=3cm,right=3cm,marginparwidth=1.75cm]{geometry}
\usepackage{amsmath, latexsym, amssymb}
\usepackage{amsthm}
\newtheorem{thm}{Theorem}
\newtheorem{conj}[thm]{Conjecture}
\theoremstyle{definition}
\newtheorem{exam}{Example}
\usepackage{graphicx}
\usepackage[colorinlistoftodos]{todonotes}
\usepackage[colorlinks=true, allcolors=blue]{hyperref}
\usepackage{comment}
\usepackage{tikz}
\usepackage{algorithm}
\usepackage{algorithmic}

\begin{document}
\title[Saturated equiangular lines]{Saturated configuration and new large construction of equiangular lines}
\author[Y.-C.~R.~Lin]{Yen-chi~Roger~Lin}
\address{Department of Mathematics, National Taiwan Normal University, Taiwan}
\email{yclinpa@gmail.com (Y.-C.~R.~Lin)}
\author[W.-H.~Yu]{Wei-Hsuan Yu}
\address{Department of Mathematics, National Central University, Taiwan}
\email{u690604@gmail.com (W.-H.~Yu)}
\date{\today}

\begin{abstract}
A set of lines through the origin in  Euclidean space is called equiangular when any pair of lines from the set intersects with each other at a common angle. 
We study the maximum size of equiangular lines in Euclidean space and use graph theoretic approach to prove that all the currently known construction for maximum equiangular lines in $\mathbb R^d$ cannot add another line to form a larger equiangular set of lines if $14 \leq d \leq 20$ and $d \neq 15$.
We give new constructions of large equiangular lines which are 248 equiangular lines in $\mathbb R^{42}$, 200 equiangular lines in $\mathbb{R}^{41}$, 168 equiangular lines in $\mathbb{R}^{40}$, 152 equiangular lines in $\mathbb R^{39}$ with angle $1/7$, and 56 equiangular lines in $\mathbb R^{18}$ with angle $1/5$.
\end{abstract}

\maketitle

\section{Introduction}

A set of lines through the origin in Euclidean space is called \emph{equiangular} when any pair of lines for the set intersects with each other at a common angle. The study of equiangular lines often refers to equiangular tight frames (ETFs) which are the maximum separation codes on sphere attaining Welch bounds \cite{welch1974lower,fickus2018tremain, fickus2016equiangular}. ETFs are highly related to strongly regular graphs \cite{barg2015finite,waldron2009construction}. The problem of determining the maximum cardinality $N(d)$ of equiangular lines in $\mathbb R^d$ was extensively studied for the last 70 years.It is not hard to see that $N(2) = 3$ which is realized by the diagonal of a regular hexagon. The problem becomes not trivial for the higher dimension. Haantjes~\cite{haantjes1948} in 1948 who first showed that $N(3) = N(4) = 6$. The maximum equiangular line in $\mathbb{R}^3$ is the six diagonals of icosahedron with angle $\arccos (\frac 1 {\sqrt 5})$. In 1966, van~Lint and Seidel~\cite{vanlint1966} determined $N(d)$ for $5 \leq d \leq 7$. In 1973, Lemmens and Seidel~\cite{lemmens1973} extended the knowledge up to dimension $23$. The methods to give upper bounds for equiangular lines are diverse. Barg and Yu \cite{barg2014} used the semidefinite programming method. Greaves \emph{et al.}~\cite{greaves2016} studied it by the analysis of Seidel matrix. Balla \emph{et al.}~\cite{balla2016} relied on probability and Ramsey theory to achieve asymptotic bounds. Glazyrin \emph{et al.} ~\cite{glazyrin2018upper} suggested zonal spherical function on derived sets to obtain better bounds for infinitely many dimensions and so many other methods. Gerzon~\cite{lemmens1973} proved that $N(d) \leq \frac {d(d+1)} 2$. However, most of the cases are far away from this bound. Currently only dimensions $d=2,3,7$ and $23$ attain this bound and the equality holds if and only if tight spherical 5-designs exist \cite{delsarte1977spherical}. Notice that the maximum equiangular lines in $\mathbb R^7$ and $\mathbb R^{23}$ are universal optimal codes \cite{cohn2007universally} and they are subsets of $E_8$ root system and Leech lattice respectively. $E_8$ root system was recently shown to be the solution of sphere packing in $\mathbb{R}^8$ \cite{viazovska2017sphere}, so as Leech lattice for $\mathbb R^{24}$ \cite{cohn2017sphere}.    

To the best of our knowledge, the ranges of $N(d)$ for $2 \leq d \leq 43$ are listed in Table~\ref{tb:smallnd} (see~\cite{azarija2016,barg2014,greaves2018equiangular,szollosi2017,yu2015}).

\begin{table}[h]
	\centering
    \caption{Maximum cardinalities of equiangular lines for small dimensions}
    \label{tb:smallnd}
    \begin{tabular}{c|ccccccccc}
   		$d$    & 2 & 3--4 & 5  & 6  & 7--13 & 14     & 15 & 16     & 17 \\ 
        $N(d)$ & 3 & 6    & 10 & 16 & 28    & 28--29 & 36 & 40--41 & 48--49 \\
        \hline\hline
        $d$    & 18     & 19     & 20     & 21  & 22  & 23--41 & 42       & 43 \\
        $N(d)$ & 56--60 & 72--75 & 90--95 & 126 & 176 & 276    & 276--288 & 344
    \end{tabular}
\end{table}
Notice that the lower bound of $\mathbb{R}^{18}$ is $56$, which is a new construction shown in this paper (see Section ~\ref{sec:rand-gen}).
We are interested in the open cases of $N(d)$ for all $d \leq 43$. Though we cannot determine the value of $N(d)$, but our results give the hint of that all the current best known construction of equiangular lines should be maximal in that dimension (see Conjecture~\ref{conj:14-20}). 

A line through the origin in Euclidean space can be represented by any one of the opposite unit vectors that are parallel to the line. Whenever we talk about such a line in this article, we will sometimes refer to it as a unit vector that represents the line.

\section{Determination of equiangular lines being saturated}
\label{sec:saturated}

A set $X$ of lines through the origin in $\mathbb R^d$ is called equiangular \emph{with the angle} $\alpha$ if each pair of lines in $X$ forms the angle $\arccos \alpha$.
Given such a set $X$, is it possible to find another line $\ell$ in the span of $X$ such that $\ell$ intersects each line in $X$ at the same common angle $\alpha$?
If not, then the equiangular line set $X$ is called \emph{saturated} in $\mathbb R^d$.
The answer to this question is negative if $|X|$ has reached its known upper bound, for examples the $28$ lines in $\mathbb R^7$, consisting of all permutations of the vector $(1,1,1,1,1,1,-3,-3)$ in $\mathbb R^8$ (all of which live in the $7$-dimensional subspace $\sum_{i=1}^8 x_i = 0$ of $\mathbb R^8$).
In the cases where $N(d)$ has not been determined (for instance $d = 14$, $16$--$20$, and $42$), it is not an easy task to determine whether another line could be added to the current known construction of equiangular sets of lines.

We propose the following algorithm to answer this question in small dimensions. 
First we choose a subset $X' = \{ v_1, v_2, \dots, v_d \}$ of $X$ that forms a basis of the span of $X$.
Then we find the set $V_0$ of all unit vectors in the span of $X$ whose inner products with each of vectors in $X'$ are $\pm \alpha$.
Since we are looking for lines, we only need one vector from each pair of opposite vectors $w$ and $-w$ in $V_0$; call this set of representatives $V$.
Note that the difference set $X \setminus X'$ must be a subset of $V$.
Now we construct a simple graph $G$ whose vertex set is $V$, and two vertices $v$ and $v'$ in $G$ are adjacent if and only if $\langle v, v' \rangle = \pm \alpha$.
Then a saturated equiangular lines that contains $X'$ has the cardinality $d + \omega(G)$, where $\omega(G)$ is the clique number of $G$.
Although this does not directly answer the question whether or not $X$ is saturated in that dimension, the number $d + \omega(G)$ is still an upper bound for the number of equiangular lines that contains $X$.
If $d + \omega(G) = |X|$, then we can conclude that $X$ is saturated.

\begin{algorithm}[H]
\begin{algorithmic}[1]
    \STATE Find a basis $B = \{ b_i \colon 1 \leq i \leq r \}$ in $X$ that spans $\mathbb R^r$
        
    \STATE Solves for all unit vectors $v_i$ in $\mathbb R^r$ which intersect with every vectors in $B$ with angle $\pm\alpha$.
        
                
    \STATE Make a graph $G = (V,E)$, $V = \{ v_i \}$; $(v_i, v_j) \in E$ if and only if $\langle v_i, v_j \rangle = \pm \alpha$
        
    \STATE Compute $N = |B| + \omega(G)$ (the clique number of $G$)
    
    \IF{$|X| = N$}
        \STATE return $X$ is saturated.
    \ENDIF
    
\end{algorithmic}
\end{algorithm}

Using this algorithm, we are able to establish the following result, which does not seem to appear anywhere in the literature.

\begin{thm}
  \label{thm:7-max}
  The following sets are saturated equiangular lines:
  \begin{itemize}
  \item The $28$ lines in $\mathbb R^{14}$, in Tremain~\cite{tremain2008}, p.~24 (or below).
  \item The $40$ lines in $\mathbb R^{16}$, in Tremain~\cite{tremain2008}, p.~25.
  \item The $48$ lines in $\mathbb R^{17}$, in Lemmens and Seidel~\cite{lemmens1973}, section~2.
  \item The $54$ lines in $\mathbb R^{18}$, in Sz\"oll\H osi~\cite{szollosi2017}.
  \item The $72$ lines in $\mathbb R^{19}$, constructed by Asche (see Taylor~\cite{taylor1971}, Theorem~8.1, or below).
  \item The $90$ lines in $\mathbb R^{20}$, constructed by Taylor~\cite{taylor1971}, Theorem~8.2 (or below).
  \end{itemize}
\end{thm}

\begin{exam}[from~\cite{tremain2008}, p.~24]
	A construction of $28$ equiangular lines in $\mathbb R^{14}$ using $(7,3,1)$-designs was given by Tremain.
	
	\begin{figure}[h]
	\centering
    \begin{tikzpicture}[scale=0.5]
      \draw[gray!30!white] (0,0) grid (28,14);
      \foreach \x in {0,1,...,6}
        \foreach \y in {0,1,2,3}
          \node at (28-4*\x-\y-0.5, \x+0.5) {$\star$};
      \foreach \x in {1,2,17,20,25,27}
        \node at (\x-0.5,14-0.5) {$\circ$};
      \foreach \x in {3,4,18,19,26,28}
        \node at (\x-0.5,14-0.5) {$\bullet$};
      \foreach \x in {1,3,5,6,21,24}
        \node at (\x-0.5,13-0.5) {$\circ$};
      \foreach \x in {2,4,7,8,22,23}
        \node at (\x-0.5,13-0.5) {$\bullet$};
      \foreach \x in {5,7,9,10,25,28}
        \node at (\x-0.5,12-0.5) {$\circ$};
      \foreach \x in {6,8,11,12,26,27}
        \node at (\x-0.5,12-0.5) {$\bullet$};
      \foreach \x in {1,4,9,11,13,14}
        \node at (\x-0.5,11-0.5) {$\circ$};
      \foreach \x in {2,3,10,12,15,16}
        \node at (\x-0.5,11-0.5) {$\bullet$};
      \foreach \x in {5,8,13,15,17,18}
        \node at (\x-0.5,10-0.5) {$\circ$};
      \foreach \x in {6,7,14,16,19,20}
        \node at (\x-0.5,10-0.5) {$\bullet$};
      \foreach \x in {9,12,17,19,21,22}
        \node at (\x-0.5, 9-0.5) {$\circ$};
      \foreach \x in {10,11,18,20,23,24}
        \node at (\x-0.5, 9-0.5) {$\bullet$};
      \foreach \x in {13,16,21,23,25,26}
        \node at (\x-0.5, 8-0.5) {$\circ$};
      \foreach \x in {14,15,22,24,27,28}
        \node at (\x-0.5, 8-0.5) {$\bullet$};
      \foreach \x in {1,2,...,28}
        \node at (\x-0.5,14.7) [scale=0.75] {$\x$};
      \foreach \x in {1,2,...,14}
        \node at (-0.7,14.5-\x) [scale=0.75] {$\x$};
      \node at (14,-1.5) {$\circ = \sqrt{1/5}, \quad \bullet = -\sqrt{1/5}, \quad \star = \sqrt{2/5}$};
    \end{tikzpicture}
    \caption{Graphic representation of 28 equiangular lines in $\mathbb R^{14}$}
    \label{fig:14-28}
	\end{figure}

Figure~\ref{fig:14-28} is a graphic representation on 28 column vectors in $\mathbb R^{14}$: each $\circ$, $\bullet$, and $\star$ shall be replaced by $\sqrt{1/5}$, $-\sqrt{1/5}$, and $\sqrt{2/5}$, respectively; empty squares shall be filled with $0$.
One checks immediately that this indeed gives $28$ equiangular lines in $\mathbb R^{14}$ with angle $1/5$.
The following computation is executed in {\tt Sage}~\cite{sagemath}.
We label these $28$ vectors by $w_1, w_2, \dots, w_{28}$ from left to right.
The vectors $\{ w_{2k} \colon k = 1, 2, \dots, 14 \}$ form a basis for $\mathbb R^{14}$.
Let $C$ be the following set of vectors in $\mathbb R^{14}$:
\begin{equation*}
	C := \bigl\{ v \in \mathbb R^{14} \colon \langle v, w_{2k} \rangle = \pm \frac15, \,\, \forall\,\, 1\leq k \leq 14; \langle v, w_2 \rangle = \frac15 \bigr\}.
\end{equation*}
Among those $2^{13}$ vectors in $C$, there are $378$ unit vectors; call the collection of these unit vectors $V$.
Using these vectors in $V$ as vertices, a simple graph $G$ is constructed by connecting $v$, $v' \in G$ whose inner product is $\pm \frac15$. Finally we verify that the clique number $\omega(G) = 14$, which means that the maximum cardinality of equiangular lines in $\mathbb R^{14}$ that contains $X'$ is $|X'| + \omega(G) = 28$.
This implies that the $28$ equiangular lines above in $\mathbb R^{14}$ is saturated.
This computation in effect reduces the number of combinations of $\{ \frac15, -\frac15 \}$-inner products to be checked from $2^{27}$ to $2^{13}$.
\end{exam}

\begin{exam}[\cite{taylor1971}, Theorems~8.1 and 8.2]
	We hereby give the construction of $90$ equiangular lines in $\mathbb R^{20}$ and $72$ equiangular lines in $\mathbb R^{19}$ by Asche, and verify that these sets are saturated.

	The constructions come from the Witt design (also known as the Steiner triple system $S(24,8,5)$, see~\cite{witt1937}).
    We first list all the $8$-subsets of $\{1, 2, \dots, 24 \}$ in lexicographical order,
    and any such subset which differs from the some subset already found in fewer than $4$ elements is discarded.
    This procedure picks out $759$ $8$-subsets of $\{ 1, 2, \dots, 24 \}$, with the first subset being $\{1, 2, 3, 4, 5, 6, 7, 8 \}$;
    let $\mathcal C$ be the collection of these $759$ subsets. 
    Let $\{ e_i \colon i = 1, \dots, 24 \}$ denote the standard basis for $\mathbb R^{24}$,
    and let $e_\Sigma := \sum_{i=1}^{24} e_i$.
    For any $E \in \mathcal C$ with $1 \in E$, define $f(E) := \bigl(4 \sum_{i\in E} e_i - 4 e_1 - e_\Sigma\bigr)/\sqrt{80}$.
	Define
    \begin{align}
        c &:= 4 e_1 + e_{\Sigma}; \notag \\
    	c_1 &:= e_2 + e_3 + e_{10} + e_{12} + e_{13} + e_{14} + e_{21} + e_{24},;
        \label{eq:vector-c} \\
        c_2 &:= e_2 + e_3 + e_6 + e_7 + e_{18} + e_{19} + e_{22} + e_{23}.  \notag     
    \end{align}
    There are $90$ sets $E_1, \dots, E_{90}$ in $\mathcal C$ such that $1 \in E_i$ and $f(E_i)$ is perpendicular to each of $e_1 - e_2$, $c$, $c_1$, and $c_2$, for all $1 \leq i \leq 90$.
    Moreover, it is readily checked that the unit vectors $f(E_i)$, $1 \leq i \leq 90$, form equiangular lines with common angle $\frac15$ and live in a $20$-dimensional subspace $W$ of $\mathbb R^{24}$.
    The list of these $8$-sets $E_i$ is given in Table~\ref{tb:20-90} by the lexicographic order.
    
    \begin{table}
    \centering
    \caption{The 90 subsets of $\{1,2,\dots,24\}$ used to construct equiangular lines in $\mathbb R^{20}$}
    \label{tb:20-90}
    \begin{tabular}{|l|l|l|} \hline
	1, 3, 4, 5, 9, 15, 18, 24 & 1, 4, 7, 8, 13, 15, 22, 24 & 1, 6, 7, 8, 10, 15, 20, 21 \\ 
 	1, 3, 4, 5, 10, 16, 17, 23 & 1, 4, 7, 8, 14, 16, 21, 23 & 1, 6, 7, 8, 11, 14, 17, 24 \\ 
	1, 3, 4, 5, 11, 13, 20, 22 & 1, 4, 9, 10, 14, 16, 18, 19 & 1, 6, 9, 10, 11, 14, 20, 22 \\ 
 	1, 3, 4, 6, 9, 16, 20, 21 & 1, 4, 9, 12, 17, 18, 21, 23 & 1, 6, 9, 11, 12, 16, 18, 24 \\ 
 	1, 3, 4, 7, 11, 15, 17, 21 & 1, 4, 9, 12, 19, 20, 22, 24 & 1, 6, 9, 13, 14, 16, 17, 23 \\ 
 	1, 3, 4, 8, 9, 14, 17, 22 & 1, 4, 10, 11, 17, 18, 22, 24 & 1, 6, 9, 15, 17, 21, 22, 24 \\ 
 	1, 3, 4, 8, 11, 16, 19, 24 & 1, 4, 10, 11, 19, 20, 21, 23 & 1, 6, 10, 11, 12, 15, 17, 23 \\ 
 	1, 3, 4, 8, 12, 15, 20, 23 & 1, 4, 11, 12, 13, 15, 18, 19 & 1, 6, 10, 16, 17, 19, 20, 24 \\ 
 	1, 3, 5, 6, 14, 15, 17, 20 & 1, 4, 13, 16, 17, 19, 21, 22 & 1, 6, 11, 13, 17, 18, 20, 21 \\ 
 	1, 3, 5, 7, 9, 11, 14, 16 & 1, 4, 13, 16, 18, 20, 23, 24 & 1, 6, 11, 14, 15, 16, 19, 21 \\ 
 	1, 3, 5, 8, 9, 10, 19, 20 & 1, 4, 14, 15, 17, 19, 23, 24 & 1, 6, 12, 13, 15, 16, 20, 22 \\ 
 	1, 3, 5, 8, 11, 12, 17, 18 & 1, 4, 14, 15, 18, 20, 21, 22 & 1, 7, 9, 10, 11, 15, 19, 24 \\ 
 	1, 3, 5, 8, 15, 16, 21, 22 & 1, 5, 6, 7, 9, 13, 20, 24 & 1, 7, 9, 10, 12, 16, 20, 23 \\ 
 	1, 3, 6, 8, 9, 11, 13, 15 & 1, 5, 6, 7, 12, 16, 17, 21 & 1, 7, 9, 11, 12, 13, 17, 22 \\ 
 	1, 3, 7, 8, 13, 16, 17, 20 & 1, 5, 6, 8, 9, 14, 18, 21 & 1, 7, 9, 13, 15, 16, 18, 21 \\ 
 	1, 3, 9, 11, 17, 20, 23, 24 & 1, 5, 6, 8, 10, 13, 17, 22 & 1, 7, 9, 14, 17, 19, 20, 21 \\ 
 	1, 3, 9, 12, 15, 16, 17, 19 & 1, 5, 6, 8, 12, 15, 19, 24 & 1, 7, 10, 14, 15, 16, 17, 22 \\ 
 	1, 3, 10, 11, 15, 16, 18, 20 & 1, 5, 7, 8, 10, 16, 18, 24 & 1, 7, 11, 13, 14, 15, 20, 23 \\ 
 	1, 4, 5, 6, 10, 12, 18, 20 & 1, 5, 7, 8, 11, 13, 19, 21 & 1, 7, 11, 16, 20, 21, 22, 24 \\ 
 	1, 4, 5, 6, 13, 15, 21, 23 & 1, 5, 7, 8, 12, 14, 20, 22 & 1, 7, 12, 15, 17, 18, 20, 24 \\ 
 	1, 4, 5, 6, 14, 16, 22, 24 & 1, 5, 9, 10, 11, 13, 18, 23 & 1, 8, 9, 10, 12, 15, 18, 22 \\ 
 	1, 4, 5, 7, 9, 10, 21, 22 & 1, 5, 9, 13, 14, 15, 19, 22 & 1, 8, 9, 11, 12, 14, 19, 23 \\ 
 	1, 4, 5, 7, 11, 12, 23, 24 & 1, 5, 9, 16, 19, 21, 23, 24 & 1, 8, 9, 13, 17, 18, 19, 24 \\ 
 	1, 4, 5, 7, 13, 14, 17, 18 & 1, 5, 10, 11, 12, 16, 19, 22 & 1, 8, 9, 13, 20, 21, 22, 23 \\ 
 	1, 4, 6, 7, 9, 12, 14, 15 & 1, 5, 10, 15, 17, 18, 19, 21 & 1, 8, 10, 13, 15, 16, 19, 23 \\ 
 	1, 4, 6, 7, 10, 11, 13, 16 & 1, 5, 10, 15, 20, 22, 23, 24 & 1, 8, 10, 14, 17, 18, 20, 23 \\ 
 	1, 4, 6, 8, 9, 10, 23, 24 & 1, 5, 11, 14, 17, 21, 22, 23 & 1, 8, 11, 13, 14, 16, 18, 22 \\ 
 	1, 4, 6, 8, 11, 12, 21, 22 & 1, 5, 11, 14, 18, 19, 20, 24 & 1, 8, 11, 15, 18, 21, 23, 24 \\ 
 	1, 4, 6, 8, 13, 14, 19, 20 & 1, 5, 12, 13, 17, 19, 20, 23 & 1, 8, 12, 16, 17, 22, 23, 24 \\ 
 	1, 4, 7, 8, 10, 12, 17, 19 & 1, 5, 12, 14, 15, 16, 18, 23 & 1, 8, 12, 16, 18, 19, 20, 21 \\ \hline
    \end{tabular}
    \end{table}
    
    It can be checked that $X' := \{ f(E_j) \colon j \in J \}$ forms a basis for $W$, where $J = \{ 6, 7, 13, 19, 21, 24,$ $27, 34, 43, 45, 48, 52, 57, 61, 66, 70, 74, 80, 82, 89 \}$.
    Let $C$ be the following set of vectors in $W$:
	\begin{equation*}
	C := \bigl\{ v \in W \colon \langle v, f(E_j) \rangle = \pm \frac15, \,\, \forall\,\, j \in J; \langle v, f(E_6) \rangle = \frac15 \bigr\}.
\end{equation*}
Among those $2^{19}$ vectors in $C$, there are only $70$ vectors of unit length; in fact they are the remaining $70$ unit vectors $f(E_k)$ (or their opposites) for $k \in \{1, 2, \dots, 90 \} \setminus J$.
From here we can conclude that these $90$ lines form a saturated equiangular line set in $\mathbb R^{20}$. Note that our procedure again reduces the number of vectors to be checked from $2^{89}$ to $2^{19}$.

Inside the above $90$ lines, we can pick out $72$ lines by discarding those $18$ $8$-sets that contains $3$ from Table~\ref{tb:20-90}.
The resulting vectors $f(E_i)$, $i = 19, \dots, 90$, are also perpendicular to $e_1 - e_3$ in $\mathbb R^{24}$.
Therefore those $72$ equiangular lines live in a $19$-dimensional subspace of $\mathbb R^{24}$.

\end{exam}

\section{Construction of a large equiangular subset of lower rank}
\label{sec:rand-gen}

A large equiangular set can be found from a larger set from higher dimensional spaces.
For example, $48$ lines in $\mathbb R^{17}$, $54$ lines in $\mathbb R^{18}$, $72$ lines in $\mathbb R^{19}$, $90$ lines in $\mathbb R^{20}$, $126$ lines in $\mathbb R^{21}$, and $176$ lines in $\mathbb R^{22}$ can all be found among the $276$ equiangular lines in $\mathbb R^{23}$, sitting inside various lower dimensional subspaces. Similar stories also happen to $\mathbb{R}^7$. The $28$ equiangular lines in $\mathbb{R}^7$ contain maximum size of equiangular lines in $\mathbb{R}^6$ (16 lines) and $\mathbb{R}^5$ (10 lines).
A linearly independent subset $T$ of vectors inside an equiangular set $X$ generates the maximal subset of $X$ that contains in the span of $T$, usually of lower rank.
In this section we mention two such constructions.

\subsection{248 equiangular lines in \texorpdfstring{$\mathbb R^{42}$}{R42} with angle \texorpdfstring{$1/7$}{1/7}}
\label{sec:42-248}

The existence of $344$ equiangular lines in $\mathbb R^{43}$ with angle $1/7$ follows from Taylor's result on the doubly transitive group $\text{P$\Gamma$U}(3, 7^2)$ \cite{taylor1971}, and can also be constructed from the strongly regular graph $\operatorname{SRG}(344, 168, 92, 72)$ \cite{cameron2004strongly}, which induces the Gram matrix $G = [\langle v_i, v_j \rangle]_{i,j=1}^{344}$ of these lines.
We \emph{randomly} select $42$ columns from $G$ and verify that they form a linearly independent set of vectors.
Then we collect the column vectors of $G$ which belong to the span of these vectors; call this collection $X$.
By picking out the corresponding rows and columns of $X$ from $G$, the resulting matrix is the Gram matrix of equiangular unit vectors of rank $42$ and angle $1/7$.
The best result among a few thousand runs of this experiment gave $248$ equiangular lines in $\mathbb R^{42}$ with angle $1/7$.

For sake of comparison, we recall the inequality~(\ref{eq:relative-bound}), which is the so-called \emph{relative bound} for equiangular lines.
\begin{thm}[\cite{vanlint1966}, p.342]
\label{thm:relative-bound}
Let $X$ be an equiangular set with angle $\alpha$ in $\mathbb R^r$.
If $r < \frac{1}{\alpha^2}$, then
\begin{equation}
\label{eq:relative-bound}
|X| \leq R(r,\alpha) := \frac{r (1 - \alpha^2)}{1 - r \alpha^2}.
\end{equation}
\end{thm}

Proceeding in a similar fashion, we look for large subsets in $\mathbb R^d$ of the $344$ equiangular lines in $\mathbb R^{43}$.  
The best results\footnote{Examples on these Gram matrices can be downloaded at \url{http://math.ntnu.edu.tw/~yclin/Gram1-7/}} are listed in Table~\ref{tab:4x-1/7}, with a comparison with the relative bounds (which are also the best upper bounds so far).
Although our constructions do not reach the relative bounds, we doubt if there are more equiangular lines with angle $1/7$ in these dimensions. Notice that if our construction are the maximum in that dimension and angle, then there will be new results for nonexistence of two associated strongly regular graphs (with 288, 246 vertices, respectively).

\begin{table}
    \centering
    \caption{Relative bounds and the sizes of known constructions of equiangular set with angle $1/7$ in $\mathbb R^d$}
    \label{tab:4x-1/7}
    \begin{tabular}{c|cccc}
        $d$        &  42 &  41 &  40 &  39 \\ \hline
        $R(d,1/7)$ & 288 & 246 & 213 & 187 \\ 
        Found      & 248 & 200 & 168 & 152 
    \end{tabular}
\end{table}

\subsection{56 equiangular lines in \texorpdfstring{$\mathbb R^{18}$}{R18} with angle \texorpdfstring{$1/5$}{1/5}}
\label{sec:18-56}

Among the $72$ equiangular lines in $\mathbb R^{19}$ with angle $1/5$ given in Example~2 of Section~\ref{sec:saturated}, we \emph{randomly} select $18$ of them and collect all the vectors from those $72$ vectors that fall into the span of these $18$ vectors.
After a few hundred runs of the experiments, the best result we find is a collection of $56$ equiangular lines of rank $18$ with angle $1/5$.
Specifically, there are two kinds of such $56$ vectors.
Let
\begin{align}
    u_1 &:= e_4 + e_5 + 6e_6 - 3e_7 + e_8 + e_9 + 6 e_{10} + e_{11} - 3 e_{12} + 6 e_{13} - 3 e_{14} + e_{15} + e_{16} + e_{17} \notag \\
    &\phantom{:= }   - 3 e_{18} - 3 e_{19} - 8 e_{20} - 3 e_{21} - 3 e_{22} + 6 e_{23} - 3 e_{24};  \label{eq:vector-u} \\
    u_2 &:= 5 e_4 + 5 e_5 - 3 e_6 - 3 e_7 - 4 e_8 + 5 e_9 - 4 e_{11} - 4 e_{15} + 5 e_{16} - 4 e_{17} - 3 e_{18} + 6 e_{19} - 4 e_{20}  \notag \\
    &\phantom{:= }  + 6 e_{22} - 3 e_{23}. \notag
\end{align}
The collections we find are perpendicular to either $\{ c, c_1, c_2, e_1 - e_2, e_1 - e_3, u_1 \}$ or $\{ c, c_1, c_2, e_1 - e_2, e_1 - e_3, u_2 \}$ (the vectors $c$, $c_1$, $c_2$ are listed in (\ref{eq:vector-c}); the coordinates of $u_1$ and $u_2$ in (\ref{eq:vector-u}) may be permuted in specific ways.)
This finding raises the lower bound of $N(18)$ from $54$ to $56$, see Table~\ref{tb:smallnd}.
Also from the method described in Section~\ref{sec:saturated}, we confirm that both configurations of $56$ equiangular lines in $\mathbb R^{18}$ are saturated.

\section{Discussions}
All other equiangular sets of lines in Theorem~\ref{thm:7-max} are checked to be saturated using the similar procedures to the Examples above.
But it took too long when we tried to examine $248$ equiangular lines in $\mathbb R^{42}$.
Our algorithm would require to pick out unit vectors among $2^{41}$ (roughly $2.2$ trillion) possible candidates.

In light of Theorem~\ref{thm:7-max}, we propose this following conjecture.

\begin{conj}
    \label{conj:14-20}
	The following table gives on the maximal equiangular lines on the specified dimensions:
    \begin{equation*}
    	\begin{array}{c|cccccc}
        	d    & 14 & 16 & 17 & 18 & 19 & 20 \\ \hline
            N(d) & 28 & 40 & 48 & 56 & 72 & 90
        \end{array}
    \end{equation*}
\end{conj}

Conjecture~\ref{conj:14-20} is coherent to Peter Casazza's conjecture which states that all the maximal sizes of equiangular lines are even numbers except for $\mathbb{R}^2$. We know that the crucial steps are to do the classification of equiangular lines. For instance, the 36 equiangular lines in $\mathbb{R}^{15}$ with angle $\frac 1 5$ has been classified to be 227 different classes \cite{spence1995regular}. We believe that there should be less than 227 different classes of 28 equingular lines in $\mathbb{R}^{14}$. If there are not so many different classes, then in conjunction of our methods, we might be able to prove $N(14)=28$.  

\subsection*{Acknowledgements}
The first author is partially supported by 107-2115-M-003-001 from Ministry of Science and Technology, Taiwan.
Part of this work is based upon work supported by the National Science Foundation under
Grant No.~DMS-1439786 while the second author was in residence at the Institute for 
Computational and Experimental Research in Mathematics in Providence, RI, during the Point configurations in Geometry, Physics and Computer Science Program.
Part of this work was done when the second author visited National Center for Theoretical Sciences (NCTS), Taiwan, in the summer of 2017.
The authors are grateful to the support of NCTS.
The authors wish to thank Prof.~E.~Bannai and Gary Greaves for their useful comments.

\bibliographystyle{amsplain}
\bibliography{equiangular}

\providecommand{\bysame}{\leavevmode\hbox to3em{\hrulefill}\thinspace}
\providecommand{\MR}{\relax\ifhmode\unskip\space\fi MR }
\providecommand{\MRhref}[2]{%
  \href{http://www.ams.org/mathscinet-getitem?mr=#1}{#2}
}
\providecommand{\href}[2]{#2}
\begin{thebibliography}{10}

\bibitem{azarija2016}
Jernej Azarija and Tilen Marc, \emph{There is no (95, 40, 12, 20) strongly
  regular graph}, arXiv preprint arXiv:1603.02032 (2016).

\bibitem{balla2016}
Igor Balla, Felix Dr{\"a}xler, Peter Keevash, and Benny Sudakov,
  \emph{Equiangular lines and spherical codes in {E}uclidean space},
  Inventiones mathematicae (2017), {\tt
  https://doi.org/10.1007/s00222-017-0746-0}.

\bibitem{barg2015finite}
Alexander Barg, Alexey Glazyrin, Kasso~A Okoudjou, and Wei-Hsuan Yu,
  \emph{Finite two-distance tight frames}, Linear Algebra and its Applications
  \textbf{475} (2015), 163--175.

\bibitem{barg2014}
Alexander Barg and Wei-Hsuan Yu, \emph{New bounds for equiangular lines},
  Contemporary Mathematics \textbf{625} (2014), 111--121.

\bibitem{cameron2004strongly}
Peter~J Cameron, \emph{Strongly regular graphs}, Topics in Algebraic Graph
  Theory \textbf{102} (2004), 203--221.

\bibitem{cohn2007universally}
Henry Cohn and Abhinav Kumar, \emph{Universally optimal distribution of points
  on spheres}, Journal of the American Mathematical Society \textbf{20} (2007),
  no.~1, 99--148.

\bibitem{cohn2017sphere}
Henry Cohn, Abhinav Kumar, Stephen~D Miller, Danylo Radchenko, and Maryna
  Viazovska, \emph{The sphere packing problem in dimension 24}, Annals of
  Mathematics (2017), 1017--1033.

\bibitem{delsarte1977spherical}
P~Delsarte, \emph{Spherical codes and designs}, Geom. Dedicata \textbf{6}
  (1977), 363--388.

\bibitem{fickus2018tremain}
Matthew Fickus, John Jasper, Dustin~G Mixon, and Jesse Peterson, \emph{Tremain
  equiangular tight frames}, Journal of Combinatorial Theory, Series A
  \textbf{153} (2018), 54--66.

\bibitem{fickus2016equiangular}
Matthew Fickus, Dustin~G Mixon, and John Jasper, \emph{Equiangular tight frames
  from hyperovals}, IEEE Transactions on Information Theory \textbf{62} (2016),
  no.~9, 5225--5236.

\bibitem{glazyrin2018upper}
Alexey Glazyrin and Wei-Hsuan Yu, \emph{Upper bounds for $s$-distance sets and
  equiangular lines}, Advances in Mathematics \textbf{330} (2018), 810--833.

\bibitem{greaves2016}
Gary Greaves, Jacobus~H. Koolen, Akihiro Munemasa, and Ferenc
  Sz{\"o}ll{\H{o}}si, \emph{Equiangular lines in {E}uclidean spaces}, Journal
  of Combinatorial Theory, Series A \textbf{138} (2016), 208--235.

\bibitem{greaves2018equiangular}
Gary~RW Greaves, \emph{Equiangular line systems and switching classes
  containing regular graphs}, Linear Algebra and its Applications \textbf{536}
  (2018), 31--51.

\bibitem{haantjes1948}
J.~Haantjes, \emph{Equilateral point-sets in elliptic two- and
  three-dimensional spaces}, Nieuw Arch. Wisk \textbf{22} (1948), no.~2,
  355--362.

\bibitem{lemmens1973}
Petrus W.~H. Lemmens and Johan~J. Seidel, \emph{Equiangular lines}, Journal of
  Algebra \textbf{24} (1973), no.~3, 494--512.

\bibitem{sagemath}
The {Sage Developers}, \emph{{S}agemath, the {S}age {M}athematics {S}oftware
  {S}ystem ({V}ersion 8.3)}, 2018, {\tt http://www.sagemath.org}.

\bibitem{spence1995regular}
Edward Spence, \emph{Regular two-graphs on 36 vertices}, Linear algebra and its
  applications \textbf{226} (1995), 459--497.

\bibitem{szollosi2017}
Ferenc Sz{\"o}ll{\H{o}}si, \emph{A remark on a construction of {D.~S.~Asche}},
  arXiv preprint arXiv:1703.04505 (2017).

\bibitem{taylor1971}
Donald~E. Taylor, \emph{Some topics in the theory of finite groups}, Ph.D.
  thesis, University of Oxford, 1971.

\bibitem{tremain2008}
Janet~C. Tremain, \emph{Concrete constructions of real equiangular line sets},
  arXiv preprint arXiv:0811.2779 (2008), 1--39.

\bibitem{vanlint1966}
Jacobus~H. van Lint and Johan~J. Seidel, \emph{Equilateral point sets in
  elliptic geometry}, Indag. Math \textbf{28} (1966), no.~3, 335--348.

\bibitem{viazovska2017sphere}
Maryna~S Viazovska, \emph{The sphere packing problem in dimension 8}, Annals of
  Mathematics (2017), 991--1015.

\bibitem{waldron2009construction}
Shayne Waldron, \emph{On the construction of equiangular frames from graphs},
  Linear Algebra and its Applications \textbf{431} (2009), no.~11, 2228--2242.

\bibitem{welch1974lower}
Lloyd Welch, \emph{Lower bounds on the maximum cross correlation of signals
  (corresp.)}, IEEE Transactions on Information theory \textbf{20} (1974),
  no.~3, 397--399.

\bibitem{witt1937}
Ernst Witt, \emph{Die 5-fach transitiven {G}ruppen von {M}athieu}, Abhandlungen
  aus dem Mathematischen Seminar der Universit{\"a}t Hamburg, vol.~12,
  Springer, 1937, pp.~256--264.

\bibitem{yu2015}
Wei-Hsuan Yu, \emph{There are no 76 equiangular lines in {R}$^{19}$}, arXiv
  preprint arXiv:1511.08569 (2015), 1--12.

\end{thebibliography}


\begin{thebibliography}{10}

\bibitem{azarija2016}
Jernej Azarija and Tilen Marc, \emph{There is no $(95, 40, 12, 20)$ strongly
  regular graph}, arXiv preprint arXiv:1603.02032 (2016).

\bibitem{barg2014}
Alexander Barg and Wei-Hsuan Yu, \emph{New bounds for equiangular lines},
  Contemporary Mathematics \textbf{625} (2014), 111--121.

\bibitem{cameron2004}
Peter~J.~Cameron, \emph{Strongly regular graphs},
  Topics in Algebraic Graph Theory
  \textbf{102} (2004), 203--221.

\bibitem{greaves2016-1}
Gary R.~W. Greaves, \emph{Equiangular line systems and switching classes
  containing regular graphs}, arXiv preprint arXiv:1612.03644 (2016).

\bibitem{haantjes1948}
J.~Haantjes, \emph{Equilateral point-sets in elliptic two- and
  three-dimensional spaces}, Nieuw Arch. Wisk \textbf{22} (1948), no.~2,
  355--362.

\bibitem{lemmens1973}
Petrus W.~H. Lemmens and Johan~J. Seidel, \emph{Equiangular lines}, Journal of
  Algebra \textbf{24} (1973), no.~3, 494--512.

\bibitem{sagemath}
The {Sage Developers}, \emph{{S}agemath, the {S}age {M}athematics {S}oftware
  {S}ystem ({V}ersion 8.3)}, 2018, obtainable at {\tt http://www.sagemath.org}.

\bibitem{szollosi2017}
Ferenc Sz{\"o}ll{\H{o}}si, \emph{A remark on a construction of {D.~S.~Asche}},
  arXiv preprint arXiv:1703.04505 (2017).

\bibitem{taylor1971}
Donald~E. Taylor, \emph{Some topics in the theory of finite groups}, Ph.D.
  thesis, University of Oxford, 1971.

\bibitem{tremain2008}
Janet~C. Tremain, \emph{Concrete constructions of real equiangular line sets},
  arXiv preprint arXiv:0811.2779 (2008).

\bibitem{vanlint1966}
Jacobus~H. van Lint and Johan~J. Seidel, \emph{Equilateral point sets in
  elliptic geometry}, Indag. Math \textbf{28} (1966), no.~3, 335--348.

\bibitem{witt1937}
Ernst Witt, \emph{Die 5-fach transitiven {G}ruppen von {M}athieu}, Abhandlungen
  aus dem Mathematischen Seminar der Universit{\"a}t Hamburg, vol.~12,
  Springer, 1937, pp.~256--264.

\bibitem{yu2015}
Wei-Hsuan Yu, \emph{There are no $76$ equiangular lines in \texorpdfstring{$\mathbb R^{19}$}{R19}}, arXiv
  preprint arXiv:1511.08569 (2015).

\end{thebibliography}
\end{document}